\documentclass{amsart}
\usepackage{amsmath}
\usepackage{amssymb}
\usepackage{amsthm}
\usepackage[msc-links,abbrev]{amsrefs}
\usepackage{hyperref}
\usepackage[noabbrev,capitalize]{cleveref}

\DeclareMathOperator{\supp}{supp}

\DeclareMathOperator{\Ric}{Ric}

\DeclareMathOperator{\Imaginary}{Im}
\DeclareMathOperator{\Real}{Re}

\newcommand{\og}{\overline{g}}

\newcommand{\htheta}{\widehat{\theta}}

\newcommand{\lp}{\langle}
\newcommand{\rp}{\rangle}
\newcommand{\lv}{\lvert}
\newcommand{\rv}{\rvert}

\newcommand{\db}{\partial_b}
\newcommand{\dbbar}{\overline{\partial}_b}

\newcommand{\bC}{\mathbb{C}}

\def\sideremark#1{\ifvmode\leavevmode\fi\vadjust{\vbox to0pt{\vss
 \hbox to 0pt{\hskip\hsize\hskip1em
 \vbox{\hsize3cm\tiny\raggedright\pretolerance10000
 \noindent #1\hfill}\hss}\vbox to8pt{\vfil}\vss}}}

\newcommand{\suchthat}{\mathrel{}\middle|\mathrel{}}

\newcommand{\comment}[1]{}

\newtheorem{theorem}{Theorem}[section]
\newtheorem{proposition}[theorem]{Proposition}

\theoremstyle{definition}

\theoremstyle{remark}

\numberwithin{equation}{section}

\usepackage{tikz-cd}

\begin{document}

\title[Lichnerowicz--Obata theorem for $\Box_b$ in three dimensions]{The Lichnerowicz--Obata theorem for the Kohn Laplacian in three dimensions}
\author{Jeffrey S. Case}
\thanks{JSC was supported by a grant from the Simons Foundation (Grant No.\ 524601)}
\address{Department of Mathematics \\ Penn State University \\ University Park, PA 16802 \\ USA}
\email{jscase@psu.edu}
\author{Paul Yang}
\thanks{PY was supported by a grant from the National Science Foundation (Grant No.\ DMS-1607091)}
\address{Department of Mathematics \\ Princeton University \\ Princeton, NJ 08540 \\ USA}
\email{yang@math.princeton.edu}
\begin{abstract}
 We prove rigidity for the Lichnerowicz-type eigenvalue estimate for the Kohn Laplacian on strictly pseudoconvex three-manifolds with nonnegative CR Paneitz operator and positive Webster curvature.
\end{abstract}
\maketitle

\section{Introduction}
\label{sec:intro}

The classical Lichnerowicz--Obata Theorem~\cites{Obata1962,Lichnerowicz1958} states that if $(M^n,g)$ is a closed Riemannian manifold with $\Ric\geq (n-1)\kappa>0$, then the first nonzero eigenvalue $\lambda_1(-\Delta)$ of the Laplacian on functions satisfies $\lambda_1(-\Delta)\geq n\kappa$, with equality if and only if $(M^n,g)$ is isometric to the round $n$-sphere of constant sectional curvature $\kappa$.  The contribution of Lichnerowicz~\cite{Lichnerowicz1958} is to prove the estimate of $\lambda_1(-\Delta)$, while the contribution of Obata is to prove rigidity; i.e.\ the characterization of equality.

Chanillo, Chiu and Yang~\cite{ChanilloChiuYang2010} proved an analogue of the Lichnerowicz eigenvalue estimate for the Kohn Laplacian $\Box_b:=2\dbbar^\ast\dbbar$.  Specifically, they proved that any closed pseudohermitian three-manifold with nonnegative CR Paneitz operator and positive Webster curvature $R$ satisfies $\lambda_1(\Box_b)\geq\inf R$, where $\lambda_1(\Box_b)$ denotes the infimum of the positive eigenvalues of $\Box_b$.  There are two subtle differences in comparison with the Riemannian case.  First, $\ker\Box_b$ is generally infinite-dimensional, as it consists of the space of CR functions.  Second, the assumption on the CR Paneitz operator --- which is a fourth-order CR invariant operator --- cannot be removed; this can be seen by considering the Rossi spheres~\cites{ChanilloChiuYang2010,Rossi1965}.  Using an embeddability result of Kohn~\cite{Kohn1985}, it follows that any closed strictly pseudoconvex three-manifold with positive CR Yamabe constant and nonnegative CR Paneitz operator --- all CR invariant assumptions --- is globally embeddable.

Li, Son and Wang~\cite{LiSonWang2015} observed that the argument of Chanillo, Chiu and Yang~\cite{ChanilloChiuYang2010} extends to higher dimensions, where a fourth-order assumption is no longer required.  More precisely, they proved that if $(M^{2n+1},T^{1,0},\theta)$, $n\geq2$, is a closed pseudohermitian manifold for which the Ricci curvature of the Tanaka--Webster connection satisfies $\Ric\geq \kappa>0$, then $\lambda_1(\Box_b)\geq\frac{2n}{n+1}\kappa$.  They also proved rigidity: If $\lambda_1(\Box_b)=\frac{2n}{n+1}\kappa$, then $(M^{2n+1},T^{1,0},\theta)$ is isometric to the round CR $(2n+1)$-sphere with $\Ric=\kappa$.  Their proof relies crucially on the fact that $\ker B=\ker\dbbar^\ast\db^\ast B$ when $n\geq2$, where $Bf:=dd_b^cf$ is the second-order operator whose kernel coincides with the space of CR pluriharmonic functions~\cite{Lee1988}.  When $n=1$, the corresponding statement one needs is that the kernel of the CR Paneitz operator coincides with the space of CR pluriharmonic functions.  It is presently unknown whether there exists a closed CR three-manifold for which the kernel of the CR Paneitz operator is strictly larger than the space of CR pluriharmonic functions.

Takeuchi~\cite{Takeuchi2019} recently showed that on embeddable strictly pseudoconvex three-manifolds, the kernel of the CR Paneitz operator coincides with the space of CR pluriharmonic functions.  Since embeddability follows from the assumptions of positive Webster scalar curvature and nonnegative CR Paneitz operator, the obstacle encountered by Li--Son--Wang does not arise.  As such, one expects a Lichnerowicz--Obata-type theorem for the Kohn Laplacian in dimension three.  The purpose of this note is to realize this expectation:

\begin{theorem}
 \label{thm:kohn_obata}
 Let $(M^3,T^{1,0},\theta)$ be a closed pseudohermitian three-manifold with nonnegative CR Paneitz operator and Webster curvature $R\geq\kappa>0$.  Then $\lambda_1(\Box_b)\geq\kappa$, with equality if and only if $(M^3,T^{1,0},\theta)$ is isometric to the standard CR three-sphere of constant Webster curvature $\kappa$
\end{theorem}

In particular, \cref{thm:kohn_obata} and the result of Li, Son and Wang~\cite{LiSonWang2015} establish CR analogues of the Lichnerowicz--Obata Theorem for the Kohn Laplacian in all dimensions.  Note that CR analogues of the Lichnerowicz--Obata Theorem for the sub-Laplacian are also known~\cites{LiWang2013,IvanovVassilev2015}.

This note is organized as follows.  In \cref{sec:bg} we define the CR Paneitz operator and discuss some of its important properties, especially how it arises in the work of Chanillo, Chiu and Yang~\cite{ChanilloChiuYang2010}.  In \cref{sec:proof} we prove \cref{thm:kohn_obata}.  Indeed, our argument is much simpler than the Li--Son--Wang argument, in part because the Webster curvature is a scalar in dimension three.
\section{Background}
\label{sec:bg}

A \emph{CR three-manifold} is a pair $(M^3,T^{1,0})$ of an orientable three-dimensional (real) manifold $M^3$ and a rank one distribution $T^{1,0}\subset TM\otimes\bC$ such that $T^{1,0}\cap T^{0,1}=\{0\}$, where $T^{0,1}:=\overline{T^{1,0}}$.  Since $M$ is orientable and $H:=\Real T^{1,0}\subset TM$ is a rank two distribution which is oriented by its complex structure, there is a nonvanishing (real) one-form $\theta$ on $M$ with $\ker\theta=H$.  Such a one-form is called a \emph{contact form} for $(M^3,T^{1,0})$.  Note that if $\theta,\htheta$ are two contact forms for $(M^3,T^{1,0})$, then there is a smooth (real-valued) function $\Upsilon\in C^\infty(M)$ such that $\htheta=e^\Upsilon\theta$.

A \emph{strictly pseudoconvex three-manifold} is a CR three-manifold $(M^3,T^{1,0})$ which admits a contact form $\theta$ for which the Levi form, $L_\theta(U,V):=-i\,d\theta(U,\bar V)$, is positive definite on $T^{1,0}$.  Note that if $L_\theta$ is positive definite on $T^{1,0}$, then so too is $L_{\htheta}$ for all contact forms $\htheta$ on $(M^3,T^{1,0})$.

A \emph{pseudohermitian three-manifold} is a triple $(M^3,T^{1,0},\theta)$ of a strictly pseudoconvex three-manifold $(M^3,T^{1,0})$ and a contact form $\theta$ for $(M^3,T^{1,0})$.  The \emph{Reeb vector field} $T$ is the unique (real) vector field on $M$ such that $\theta(T)=1$ and $d\theta(T,\cdot)=0$.  There is a unique connection, the \emph{Tanaka--Webster connection}~\cites{Webster1978,Tanaka1975}, associated to such a structure, defined as follows: Let $Z_1$ be a local frame of $T^{1,0}$ and set $Z_{\bar 1}:=\overline{Z_1}$.  Then $\{\theta,\theta^1,\theta^{\bar 1}\}$, the coframe dual to $\{T,Z_1,Z_{\bar 1}\}$, is an \emph{admissible coframe}; it satisfies $d\theta=ih_{1\bar 1}\theta^1\wedge\theta^{\bar 1}$ for some positive (real-valued) function $h_{1\bar 1}$.  We also refer to $h_{1\bar1}$ as the Levi form.  Denote by $h^{1\bar1}$ the multiplicative inverse of $h_{1\bar 1}$.  The connection form $\omega_1{}^1$ is uniquely determined by the equations
\begin{align*}
 d\theta^1 & = \theta^1 \wedge \omega_1{}^1 + h^{1\bar 1}A_{\bar1\bar1}\theta \wedge \theta^{\bar 1}, \\
 dh_{1\bar 1} & = 2\Real \omega_1{}^{1}h_{1\bar 1} ,
\end{align*}
where $A_{11}$ is a complex-valued function, the \emph{pseudohermitian torsion}.  The \emph{Tanaka--Webster connection} is determined from $\nabla T:=0$ and $\nabla Z_1 := \omega_1{}^1\otimes Z_1$ by linearity and conjugation.  The structure equation for the Tanaka--Webster connection is
\[ d\omega_1{}^1 \equiv R\,\theta^1\wedge\theta^{\bar 1} \mod \theta , \]
where $R$ is a real-valued function, the \emph{Webster curvature}.

Given a pseudohermitian three-manifold $(M^3,T^{1,0},\theta)$, we use the subscripts $1$, $\bar 1$, and $0$ to denote components of a covariant tensor field with respect to a given admissible coframe $\{\theta,\theta^1,\theta^{\bar 1}\}$.  We raise subscripts using $h^{1\bar 1}$.  Subscripts following a semicolon denote components of a covariant derivative, though the semicolon will be omitted when denoting covariant derivatives of a function.  For example, the exterior derivative $df$ of a complex function $f\in C^\infty(M;\bC)$ may be written
\[ df = f_1\theta^1 + f_{\bar 1}\theta^{\bar 1} + f_0\theta . \]

A \emph{$(0,q)$-form} on a CR three-manifold $(M^3,T^{1,0})$ is a complex-valued $q$-form which annihilates $T^{1,0}$.  Denote by $\Lambda^{0,q}$ the space of $(0,q)$-forms.  When $M$ is closed and strictly pseudoconvex, define Hermitian inner products on $C^\infty(M;\bC)$ and $\Lambda^{0,1}$ by
\begin{align}
 \label{eqn:ip} \lp f,g\rp & := \int_M f\og\,\theta\wedge d\theta , \\
 \label{eqn:inner_product} \lp \rho,\eta\rp & := \int_M h^{1\bar 1}\rho_{\bar 1}\overline{\eta_{\bar 1}}\,\theta\wedge d\theta ,
\end{align}
respectively.  We define $\dbbar\colon C^\infty(M;\bC)\to\Lambda^{0,1}$ by
\begin{equation}
 \label{eqn:dbbar}
 \dbbar f := f_{\bar 1}\theta^{\bar 1} + if^1{}_1\,\theta .
\end{equation}
It follows from the transformation formula for the Tanaka--Webster connection under change of contact form~\citelist{\cite{GoverGraham2005}*{Equation~(2.7)}\cite{Lee1986}*{Lemma~5.6}} that $\dbbar$ is CR invariant (cf.\ \cite{Takeuchi2019}*{Lemma~3.1}); i.e.\ $\dbbar$ is independent of the choice of contact form.  It follows from \cref{eqn:inner_product,eqn:dbbar} that $\dbbar^\ast\colon\Lambda^{0,1}\to C^\infty(M;\bC)$, the adjoint of $\dbbar$, is given by
\begin{equation}
 \label{eqn:dbbar_adjoint}
 \dbbar^\ast\rho = -\rho_{\bar 1;}{}^{\bar 1} .
\end{equation}
The \emph{Kohn Laplacian} (on functions), $\Box_b\colon C^\infty(M;\bC)\to C^\infty(M;\bC)$, is defined by $\Box_bf := 2\dbbar^\ast\dbbar f$.

A (complex-valued) function $f\in C^\infty(M;\bC)$ is \emph{CR} if $\dbbar f=0$.  A (real-valued) function $u\in C^\infty(M)$ is \emph{CR pluriharmonic} if locally $u=\Real f$ for some CR function $f$.  There is an alternative characterization of CR pluriharmonic functions as the kernel of a third-order $(0,2)$-form-valued differential operator: Define $P_{\bar 1}\colon C^\infty(M)\to\Lambda^{0,2}$ by
\begin{equation}
 \label{eqn:P1}
 P_{\bar 1}u := (u_1{}^1{}_{\bar 1} - iA_{\bar 1}{}^1u_1)\,\theta\wedge\theta^{\bar 1} .
\end{equation}
One easily checks that $P_{\bar1}u$ is the $(0,2)$-part of $dd_b^c u$, where $d_b^c:=\frac{i}{2}(\dbbar - \db)$, and hence $P_{\bar1}$ is CR invariant.  From this Lee concluded that $\ker P_{\bar1}$ is the space of CR pluriharmonic functions~\cite{Lee1988}*{Proposition~3.4}.

The \emph{CR Paneitz operator} on $(M^3,T^{1,0})$ is the operator $P\colon C^\infty(M)\to C^\infty(M)$ defined by $Pu\,\theta\wedge d\theta = -i\,d(P_{\bar 1}u)$.  It follows from the previous paragraph that the CR Paneitz operator is a real-valued CR invariant fourth-order differential operator.  It is immediate from \cref{eqn:P1} and the definition of $P$ that the CR Paneitz operator is formally self-adjoint.

Note that we may equivalently define $P:=-\dbbar^\ast P_{\bar 1}$, where $P_{\bar 1}$ is the $(0,1)$-form-valued differential operator $P_{\bar 1}u:=(u_1{}^1{}_{\bar 1}-iA_{\bar 1}{}^1u_1)\theta^{\bar 1}$.  Also, we may extend both $P_{\bar 1}$ and $P$ to differential operators on $C^\infty(M;\bC)$ by linearity.  It is this perspective that we will take.  In particular, $P\colon C^\infty(M;\bC)\to C^\infty(M;\bC)$ is formally self-adjoint.  We say that a closed strictly pseudoconvex three-manifold $(M^3,T^{1,0})$ has \emph{nonnegative CR Paneitz operator} if
\[ \lp Pf,f\rp \geq 0 \]
for all $f\in C^\infty(M;\bC)$, where the left-hand side is interpreted using \cref{eqn:ip}.

The key observation underlying the Lichnerowicz-type theorem for the Kohn Laplacian is the following Bochner-type formula of Chanillo, Chiu and Yang~\cite{ChanilloChiuYang2010}*{Proposition~2.1}:

\begin{proposition}
 \label{prop:bochner}
 Let $(M^3,T^{1,0},\theta)$ be a pseudohermitian three-manifold.  Then
 \[ -\frac{1}{2}\Box_b\lv\dbbar f\rv^2 = \lv f^{11}\rv^2 + \frac{1}{4}\lv\Box_bf\rv^2 - \lp\dbbar\Box_bf,\dbbar f\rp - \frac{1}{2}\lp\dbbar f,\dbbar\Box_bf\rp + R\lv\dbbar f\rv^2 - \lp P_{\bar 1}f,\dbbar f\rp \]
 for all $f\in C^\infty(M;\bC)$.
\end{proposition}
\section{Rigidity of the eigenvalue estimate}
\label{sec:proof}

Our proof of \cref{thm:kohn_obata} differs from the higher-dimensional proof given by Li, Son and Wang~\cite{LiSonWang2015} in two ways.  First, as already noted in the introduction, we need a result of Takeuchi~\cite{Takeuchi2019} to conclude that the kernel of the CR Paneitz operator coincides with the complexified space of CR pluriharmonic functions.  Second, we give a simple and direct proof that if $\lambda_1(\Box_b)=\kappa$, then the Webster curvature equals $\kappa$ and the pseudohermitian torsion vanishes identically.  The classification of such pseudohermitian manifolds then implies the conclusion.

\begin{proof}[Proof of \cref{thm:kohn_obata}]
 We begin by recalling the argument of Chanillo, Chiu and Yang~\cite{ChanilloChiuYang2010}: Recall that the spectrum of the Kohn Laplacian in $(0,\infty)$ consists only of point eigenvalues~\cite{BurnsEpstein1990}*{Theorem~1.3}.  Let $\lambda>0$ be a nonzero eigenvalue of $\Box_b$ and let $f\in C^\infty(M;\bC)$ be a nontrivial function such that $\Box_bf=\lambda f$.  Set $\kappa:=\inf R>0$.  \Cref{prop:bochner} implies that
 \begin{align*}
  0 & = \lp Pf,f\rp + \int_M \left( \lv f^{11}\rv^2 - \frac{1}{2}\lv \Box_b f\rv^2 + R\lv\dbbar f\rv^2 \right)\,\theta\wedge d\theta \\
  & \geq \frac{1}{2}\int_M \lambda(\kappa-\lambda)\lv f\rv^2
 \end{align*}
 with equality if and only if $f^{11}=0$ and $\lp Pf,f\rp=0$ and $R=\kappa$ on $\supp\,\lv\dbbar f\rv^2$.  In particular, $\lambda\geq\kappa$, and hence $\lambda_1(\Box_b)\geq\kappa$.  It follows immediately that $\Box_b$ has closed range in $L^2$.  A result of Kohn~\cite{Kohn1985} then implies that $(M^3,T^{1,0})$ is embeddable.

 We now consider the case of rigidity.  That is, suppose that $\lambda_1(\Box_b)=\kappa$ and let $f\in C^\infty(M;\bC)$ be a nontrivial function such that $\Box_bf=\kappa f$.

 First observe that $f\in\ker P_1\cap\ker P_{\bar 1}$.  Indeed, from the first paragraph it holds that $\lp Pf,f\rp=0$.  Let $u$ and $v$ be the real and imaginary parts, respectively, of $f$.  Since $P$ is a formally-self adjoint real operator, $\lp Pu,u\rp = \lp Pv,v\rp = 0$.  Therefore $u$ and $v$ are CR pluriharmonic functions~\cite{Takeuchi2019}*{Proof of Theorem~1.1}.  Since $P_1$ and $P_{\bar 1}$ are linear, we conclude that $f\in\ker P_1\cap\ker P_{\bar 1}$.

 Next observe that the set $U:=\left\{ p\in M\suchthat \dbbar f(p)\not=0\right\}$ is dense in $M$.  Indeed, let $V$ be the interior of $M\setminus U$ and suppose that $V\not=\emptyset$.  Since $\Box_bf=0$ on $V$, we have that $f=0$ on $V$ as well.  Thus $\Real f$ and $\Imaginary f$ both vanish on $V$.  By the above paragraph, $\Real f$ and $\Imaginary f$ are CR pluriharmonic functions on $V$.  Since $M$ is embeddable, the unique continuation property for holomorphic functions~\cite{BaouendiEbenfeltRothschild1999}*{Chapter~VII} implies that $\Real f=0$ and $\Imaginary f=0$ on $M$, contradicting the assumption $f\not=0$.  In particular, since $R=\kappa$ on $U$, it holds that $(M^3,T^{1,0},\theta)$ has constant Webster curvature $\kappa$.

 Now, since $P_1f=0$ and $\Box_bf=\kappa f$, it holds that
 \begin{equation}
  \label{eqn:P1-consequence}
  iA_{11}f^1 = \frac{\kappa}{2}f_1 .
 \end{equation}
 Combining \cref{eqn:P1-consequence} with the fact $f^{11}=0$ yields
 \[ if^1A_{11;}{}^{11} = \frac{\kappa}{2}f_1{}^{11} . \]
 Since $P_{\bar 1}f=0$, we then deduce that
 \[ if^1A_{11;}{}^{11} = \frac{i\kappa}{2}A^{11}f_1  . \]
 Applying \cref{eqn:P1-consequence} again yields
 \[ iA_{11;}{}^{11}f^1 = -\lv A_{11}\rv^2f^1 . \]
 Since $U$ is dense, we conclude that $iA_{11;}{}^{11}=-\lv A_{11}\rv^2$ on $M$.  Integrating over $M$ yields $A_{11}=0$.

 Finally, since $(M,T^{1,0},\theta)$ has vanishing pseudohermitian torsion and constant positive scalar curvature $\kappa$, it is isometric to a quotient of the round CR three-sphere $S_\kappa^3$ with constant Webster scalar curvature $\kappa$~\cite{Tanno1969}*{Proposition~4.1}.  It is known~\cite{Folland1972} that the eigenspace $\left\{ f\in C^\infty(S_\kappa^3) \suchthat \Box_bf = \kappa f \right\}$ is spanned by $z$ and $w$, the holomorphic coordinates on $\bC^2$.  Since there is no nontrivial linear combination of $z$ and $w$ which descends to a quotient of $S_\kappa^3$, we conclude that $(M^3,T^{1,0},\theta)$ is isometric to $S_\kappa^3$.
\end{proof}

\subsection*{Acknowledgements} The authors thank Hung-Lin Chiu for helpful conversations.

\bibliography{bib}
\end{document}